\documentclass[12pt]{amsart}
\usepackage{amsmath,amssymb,latexsym,enumitem}
\usepackage{fullpage}
\usepackage{amscd}
\theoremstyle{plain}
\newtheorem{theorem}{Theorem}
\newtheorem{proposition}{Proposition}

\newtheorem{definition}{Definition}
\newtheorem{conjecture}{Conjecture}

\def\A{\mathbb{A}}
\def\C{\mathbb{C}}

\def\Hom{\text{Hom}}

\begin{document}

\title{Doubling Constructions for Covering Groups and Tensor Product $L$-Functions}
\author{Yuanqing Cai}
\author{Solomon Friedberg}
\author{David Ginzburg}
\author{Eyal Kaplan}
\address{Cai and Friedberg:  Department of Mathematics, Boston College, Chestnut Hill MA 02467-3806, USA}
\email{yuanqing.cai@bc.edu}
\email{solomon.friedberg@bc.edu}
\address{Ginzburg: School of Mathematical Sciences, Tel Aviv University, Ramat Aviv, Tel Aviv 6997801,
Israel}
\email{ginzburg@post.tau.ac.il}
\address{Kaplan: Department of Mathematics, The Ohio State University, Columbus, OH 43210, USA}
\email{kaplaney@gmail.com}

\thanks{This work was supported by the BSF,
grant number 2012019, and by the NSF, grant number 1500977 (Friedberg).}
\begin{abstract}
This is a research announcement concerning a series of constructions obtained by applying the ``doubling method"
from the theory of automorphic forms to covering groups.  Using these constructions, we obtain partial
tensor product $L$-functions attached to generalized Shimura lifts, which may be defined in a natural way
since at almost all places the representations are unramified principal series.
\end{abstract}

\maketitle

\section{Introduction}\label{intro}
In their famous paper {\sl $L$-Functions for the Classical Groups} (in \cite{G-PS-R}), Piatetski-Shapiro and
Rallis introduced a new type of global integral which
represents the standard $L$-function for any split classical group. To
describe their construction, let ${\A}$ denote the ring of adeles
of a global field $F$. Let $G$ denote any split classical group,
$\pi_1$ and $\pi_2$ denote two irreducible cuspidal automorphic
representations of $G({\A})$, and $\varphi_{\pi_1},\varphi_{\pi_2}$
be automorphic forms in the corresponding spaces. Then the global integral they introduced
is
\begin{equation}\label{global10}
\int\limits_{G(F)\times G(F)\backslash G({\A})\times G({\A})}
\varphi_{\pi_1}(g_1)\,\varphi_{\pi_2}(g_2)E((g_1,g_2),s)\,dg_1\,dg_2.
\end{equation}
Here $E(h,s)$ is an Eisenstein series defined on the adelic points of another classical
group $H$, which depends on the choice of $G$ and
contains the direct product $G\times G$ as a subgroup.
In the case $G=GL_n$, a slight modification is required to handle the issue of convergence. Since the
domain of integration is two copies of the group $G$,
this type of construction has become known as the doubling method.
Moreover,
the doubling method also works when
$G=\widetilde{Sp}_{2n}$, the two-fold metaplectic cover of the symplectic group $Sp_{2n}$.
(This was known to Piatetski-Shapiro and Rallis.)
In all cases, after unfolding the integral it is easy to
see that \eqref{global10} is not zero only if $\pi_2$ is the contragredient of $\pi_1$.

The doubling method is general in two aspects. First, as indicated
above, it is valid for all split classical groups.
Second, and maybe more important, this construction works for {\sl{any}}
irreducible cuspidal automorphic representation $\pi_1$. This is a rare
phenomenon in these types of constructions. Usually, such integrals
unfold to some special model which is afforded by some but not all cuspidal automorphic
representations, such as the Whittaker model.
(Although every cuspidal automorphic representation of $GL_n({\A})$ is globally generic, this is not the case for other classical groups.)
In the above construction, after carrying
out the unfolding process, one obtains the inner product
\begin{equation}\label{spherical0}
<\pi_1(g)\varphi_{\pi_1},\varphi_{\pi_2}>=\int\limits_{G(F)\backslash G({\A})}
\varphi_{\pi_1}(g_1g)\,\varphi_{\pi_2}(g_1)\,dg_1.
\notag
\end{equation}
This inner product is identically zero unless $\pi_2$ is the
contragredient of $\pi_1$, and then it is nonzero for some choices of data for any $\pi_1$.

In this note we describe a generalization of the doubling method to covering groups.
Let $G$ be a split classical group,
$m\ge 1$ and let $G^{(m)}({\A})$ denote a metaplectic $m$-fold cover of $G(\A)$. This
is defined provided the underlying field has a full set of $m$-th roots of unity.
When $m=1$, we understand $G^{(1)}$ to be the linear group $G$.  
Let $\pi^{(m)}$ denote
a genuine irreducible cuspidal automorphic representation of $G^{(m)}({\A})$, and let
$\tau^{(r)}$ denote a genuine irreducible cuspidal automorphic representation of
$GL_k^{(r)}({\A})$. Here $r=m$ except when $G=Sp_{2n}$
and $m$ is even, in which case we set $r=m/2$.
Via doubling, we shall describe a global construction which represents the
partial standard tensor product $L$-function attached to $\pi^{(m)}\times
\tau^{(r)}$. This partial $L$-function is the product of
local $L$-functions over all finite places of $F$ for which the representations $\pi^{(m)}$ and
$\tau^{(r)}$ are unramified. The local $L$-functions are Rankin-Selberg $L$-functions attached to the
local Shimura lifts.  These will be described in \eqref{tensor} below, and may also be obtained from
the work of Savin \cite{Sa}.

The construction we introduce is a generalization of the construction described by
the integral \eqref{global10}, keeping track of covers.  It is given by
\begin{equation}\label{global11}
\int\limits_{G(F)\times G(F)\backslash G({\A})\times G({\A})}
\varphi_1^{(m)}(g_1)\,\overline{\varphi_2^{(m)}(g_2)}\,
E^{U,\psi_U}_{\tau^{(r)}}((g_1,g_2),s)\,dg_1\,dg_2.
\end{equation}
Here, for $i=1,2$, the functions $\varphi_i^{(m)}$ are vectors in the space of $\pi^{(m)}$. There 
is an implicit choice of section $G({\A})\rightarrow G^{(m)}({\A})$ in \eqref{global11}, and the choice of $r$ in terms of $m$ implies that the integral is in fact independent of this section.
The Eisenstein series $E$ is defined on $H^{(r)}({\A})$
where $H$ is a certain suitable split classical group,
and $E^{U,\psi_U}_{\tau^{(r)}}$ is a certain Fourier coefficient of this Eisenstein
series. The induction data of $E$ depends on a
certain representation which we refer to as
a Whittaker-Speh-Shalika representation, which is defined on a cover of the general linear group.
This representation, which depends on $\tau^{(r)}$, is
defined by means of a residue of an Eisenstein series, and when $r=1$
it reduces to the well-known Speh representation, studied for example
by Jacquet \cite{J}. We discuss these representations in Section \ref{speh} below.

Unfortunately, and this is a weak point of this construction, while these
residue representations are conjectured to exist for our 
covering groups, this is not proved in general. This is because, in the absence of a general theory
of Shimura lifts for covers of the general linear group,
 the study of the necessary partial $L$-functions needed to analyze the Eisenstein series 
attached to automorphic representations on covers of the general linear group is not available.
We discuss this issue in Section \ref{res1}, and suggest some solutions in various cases. See also Conjecture \ref{conj1}.

In the sections below we sketch the general construction of the
global integrals for all classical groups $G$. We discuss the
symplectic case in more detail
and in Theorem \ref{thloc1} we state the result of the
unramified computations.
Work carrying out
the global unfolding and the computation of the local unramified integrals
for other classical groups is in progress and we shall report on it in detail in a follow-up article.

\section{Whittaker-Speh-Shalika Representations}\label{speh}
In this section we define and construct a set of representations which we will then use in the global integrals.

\subsection{Definition of Whittaker-Speh-Shalika representations}\label{def}
Let $m, a$ and $b$ denote three positive integers. Let $GL_{ab}^{(m)}({\A})$ denote an $m$-fold metaplectic cover of $GL_{ab}({\A})$, as in Kazhdan and Patterson \cite{K-P}. When $m=1$ we take $GL_{ab}^{(m)}$ to be $GL_{ab}$. Recall that the cover is split over unipotent subgroups. For example, we can choose a section 
$GL_{ab}({\A})\rightarrow GL_{ab}^{(m)}({\A})$ which is a splitting of $GL_{ab}(F)$ and 
$V({\A})$, where $V$ is the subgroup of upper triangular unipotent matrices.
Therefore notions involving unipotent orbits transfer to covering groups
{\it mutatis mutandis}.

Let $V_{a,b}$ be the unipotent radical of the standard parabolic subgroup of $GL_{ab}$ whose Levi part is $GL_b\times GL_b \times\cdots\times GL_b$. Here $GL_b$ appears $a$ times.
In terms of matrices, this group consists of all
unipotent matrices $X$ of $GL_{ab}$ of the form
\begin{equation}\label{mat1} X=
\begin{pmatrix} I&X_{1,2}&X_{1,3}&\dots&X_{1,a}\\
&I&X_{2,3}&\dots&X_{2,a}\\
&&I&\ddots&\vdots\\ &&&\ddots&X_{a-1,a}\\ &&&&I\end{pmatrix},
\qquad
X_{i,j}\in \text{Mat}_{b\times b},
\end{equation}
where $I$ is the $b\times b$ identity matrix.   Fix a nontrivial character $\psi$ of $F\backslash \A$, and
let $\psi_{a,b}$ be the character of $V_{a,b}(F)\backslash V_{a,b}(\A)$ given by
\begin{equation*}
\psi_{a,b}(X)=\psi(\text{tr}(X_{1,2}+X_{2,3}+\cdots +X_{a-1,a}))
\end{equation*}
with $X$ as in \eqref{mat1}.

Suppose $\sigma^{(m)}$ is an automorphic representation of the group $GL_{ab}^{(m)}({\A})$.
If $\varphi_{\sigma^{(m)}}$ is in the space of $\sigma^{(m)}$, then the integral
\begin{equation}\label{whspeh1}
W(\varphi_{\sigma^{(m)}})(g)=\int\limits_{V_{a,b}(F)\backslash V_{a,b}({\A})}
\varphi_{\sigma^{(m)}}(vg)\,\psi_{a,b}(v)\,dv
\end{equation}
is a Fourier coefficient corresponding to the unipotent orbit $(a^b)$. We refer to this  Fourier coefficient as
the Whittaker-Speh-Shalika coefficient of the representation  $\sigma^{(m)}$.
Notice that when $b=1$, the Fourier coefficient \eqref{whspeh1} reduces to the  well-known Whittaker coefficient, and
when it is nonzero we say that the representation $\sigma^{(m)}$ is globally generic.   For each completion $F_\nu$ of $F$
we have a similar character of $V_{a,b}(F_\nu)$ which we also denote $\psi_{a,b}$.

We make the following definition.
\begin{definition}\label{def1}
An irreducible genuine automorphic representation  $\sigma^{(m)}$ of the group $GL_{ab}^{(m)}({\A})$ is a Whittaker-Speh-Shalika representation of type $(a,b)$ if:\\
\begin{enumerate}[leftmargin=*,label=\textbf{\arabic*})]
\item\label{def:Whittaker-Speh-Shalika 1}
This representation has a nonzero Fourier coefficient corresponding to the unipotent orbit $(a^b)$, and
moreover, for all orbits which are greater than or not related to $(a^b)$, all corresponding Fourier coefficients are zero for all choices of data. In the notation of Ginzburg \cite{G1}, this statement can be written as ${\mathcal O}_{GL_{ab}^{(m)}}(\sigma^{(m)})=(a^b)$.\\
\item\label{def:Whittaker-Speh-Shalika 2}
For a finite place $\nu$, let $\sigma^{(m)}_\nu$ denote the
irreducible constituent of $\sigma^{(m)}$ at $\nu$.
Suppose that $\sigma^{(m)}_\nu$ is an unramified representation. Then ${\mathcal O}_{GL_{ab}^{(m)}}(\sigma^{(m)}_\nu)=(a^b)$.
(That is, the local analogue of part~\ref{def:Whittaker-Speh-Shalika 1} holds.) Moreover, the dimension of
\begin{equation}\label{speh1}
{\Hom}_{V_{a,b}(F_\nu)}(\sigma^{(m)}_\nu,\psi_{a,b})
\end{equation}
is one.
\end{enumerate}
\end{definition}

When $m=b=1$, it is well known that condition 
\ref{def:Whittaker-Speh-Shalika 2} in Definition \ref{def1} is satisfied. However, on covering groups not every generic representation satisfies condition \ref{def:Whittaker-Speh-Shalika 2}.

Before constructing examples of Whittaker-Speh-Shalika representations, we note that the
Fourier coefficient $W(\varphi_{\sigma^{(m)}})(g)$ enjoys an extra invariance property. Indeed, let $GL_b^{\Delta}$
denote the image of $GL_b$ inside $GL_{ab}$ under the diagonal embedding $h\mapsto h_0:=\text{diag}(h,h,\ldots,h)$.
Then $GL_b^{\Delta}$ is the stabilizer of the character
$\psi_{a,b}$ inside the group $GL_b\times GL_b \times\cdots\times GL_b$.
Consider the group
$SL_b^{\Delta}$ embedded inside $GL_b^{\Delta}$. Then the function
$W(\varphi_{\sigma^{(m)}})(g)$ is left invariant
by all matrices $h_0$ as above with $h\in SL_b({\A})$, i.e.\
$W(\varphi_{\sigma^{(m)}})(h_0g)=W(\varphi_{\sigma^{(m)}})(g)$
for all $h_0\in SL_b^{\Delta}({\A})$.
This follows since the unipotent orbit attached to $\sigma^{(m)}$ is $(a^b)$.
Moreover, if we expand along any unipotent subgroup of $SL_b({\A})$, then the nontrivial contribution to the expansion is
trivial. This follows since the nontrivial term of the expansion is
associated with a unipotent orbit which is strictly greater than $(a^b)$. See Friedberg and Ginzburg,
\cite{F-G}, Proposition 3, for details.
This means that if $m>1$, then the group $SL_b({\A})$, embedded diagonally in $GL_{ab}^{(m)}({\A})$,
must split under the $m$-fold cover. This implies $m\mid a$.

\subsection{ Construction of Whittaker-Speh-Shalika representations}
\label{construct}
In this subsection we shall construct examples of Whittaker-Speh-Shalika representations, by means of residues of Eisenstein series.

Let $k$ and $c$ be two positive integers. Denote $b=mc$.
Let $\tau^{(m)}$ denote a genuine irreducible cuspidal automorphic representation of the group $GL_k^{(m)}({\A})$,
and $\underline{s}=(s_1,\ldots,s_b)\in\C^b$. We
construct an Eisenstein series $E_{\tau^{(m)}}^{(m)}(g,\underline{s})$
defined on the group $GL_{kb}^{(m)}({\A})$ as follows. Let
$P_{k,b}$ denote the standard  parabolic subgroup of $GL_{kb}$ whose
Levi part is $GL_k\times GL_k\times\ldots\times GL_k$. As in \cite{F-G}, Section 2, we construct the induced representation
\begin{equation}\label{ind1}
\text{Ind}_{P_{k,b}^{(m)}({\A})}^{GL_{kb}^{(m)}({\A})}
(\tau^{(m)}|\cdot|^{s_1}\otimes \tau^{(m)}|\cdot|^{s_2}\otimes
\cdots\otimes \tau^{(m)}|\cdot|^{s_b})\delta_{P_{k,b}}^{1/2}.
\end{equation}
We remark that the induction process here is more complicated when $m>1$ since the $GL_k$ blocks do not commute in the covering group.
Let $E_{\tau^{(m)}}^{(m)}(g,\underline{s})$ denote the Eisenstein series associated with this induced representation.
When $m=1$, these
representations and their residues were studied by various authors. See for example \cite{J}.
When $m>1$, these Eisenstein series were constructed and studied in Suzuki \cite{Su}, Section 8.

We start with the following conjecture.
\begin{conjecture}\label{conj1}
Given  $\tau^{(m)}$ as above, the Eisenstein series $E_{\tau^{(m)}}^{(m)}(g,\underline{s})$ has a simple multi-residue at the point
$$s_1+s_2+\cdots +s_b=0;\ \ \ \ m(s_i-s_{i+1})=1;\ \ \ \ 1\le i\le b-1.$$
\end{conjecture}
We shall discuss this conjecture in Section \ref{res1} below.

Assuming Conjecture \ref{conj1}, denote by $E_{\tau^{(m)}}^{(m)}(g)$
the residue of the Eisenstein series
$E_{\tau^{(m)}}^{(m)}(g,\underline{s})$
at the above point (this depends on the choice of test vector, but we suppress this from the notation).
Let ${\mathcal L}_{\tau^{(m)}}^{(m)}$ denote the representation of $GL_{kb}^{(m)}({\A})$ generated
by all the residue functions $E_{\tau^{(m)}}^{(m)}(g)$.
Also, let  $Z_\A$ denote the center of the group $GL_{kb}^{(m)}({\A})$.
We have the following two results.
\begin{proposition}\label{prop1}
The automorphic  representation ${\mathcal L}_{\tau^{(m)}}^{(m)}$ lies in the discrete spectrum of the space $L^2(Z_\A GL_{kb}^{(m)}(F)\backslash GL_{kb}^{(m)}({\A}))$.
\end{proposition}

\begin{proposition}\label{prop2}
We have  ${\mathcal O}_{GL_{kmc}^{(m)}}
({\mathcal L}_{\tau^{(m)}}^{(m)})=((km)^c)$.
\end{proposition}
From these results we deduce that the representation
${\mathcal L}_{\tau^{(m)}}^{(m)}$ has at least one irreducible summand
which has a nonzero Fourier coefficient corresponding to the unipotent
orbit $((km)^c)$. Denote this summand by ${\mathcal E}_{\tau^{(m)}}^{(m)}$. Then we have
\begin{theorem}\label{thspeh1}
The representation ${\mathcal E}_{\tau^{(m)}}^{(m)}$ is a Whittaker-Speh-Shalika representation of type $(km,c)$.
\end{theorem}
The first condition of Definition \ref{def1} for
representation ${\mathcal E}_{\tau^{(m)}}^{(m)}$ follows from the
above discussion.
The main content of Theorem \ref{thspeh1} is
the verification of {\bf 2)} in Definition \ref{def1}.

\section{The Global Construction}\label{global}
In this section we introduce the general global integral. In the next section we treat the case of $Sp_{2n}$ in detail.

Let $m,\ n$ and $k$ denote three positive integers.
Let $G$ denote one of the split groups  $GL_n,$
$Sp_{2n},$ $SO_{2n+1}$ and $SO_{2n}$.
Let $c(n)=n$ if $G=GL_n$,
 $c(n)=2n$ if $G=Sp_{2n}$ or $G=SO_{2n}$, and $c(n)=2n+1$ if $G=SO_{2n+1}$.
Let $G^{(m)}({\A})$ denote an
$m$-fold metaplectic cover of $G({\A})$.
There is a difference  in the symplectic group case depending on the parity of $m$;
 in this section we shall assume that if $G=Sp_{2n}$ then $m$ is odd.

Depending on $G$, we introduce another classical group $H$ on which we shall construct an Eisenstein series.
Let
\begin{align*}
H=
\begin{cases}
GL_{2nkm}&G=GL_n,\\
Sp_{4nkm}&G=Sp_{2n},\\
SO_{4nkm}&G=SO_{2n},\\
SO_{2(2n+1)km}&G=SO_{2n+1}.
\end{cases}
\end{align*}
Let $P$ denote the maximal parabolic subgroup of $H$ whose Levi part is
\begin{align*}
\begin{cases}
GL_{kmc(n)}\times
GL_{kmc(n)}&G=GL_n,\\
GL_{kmc(n)}&\text{otherwise.}
\end{cases}
\end{align*}
For later use, let $U(P)$ denote the unipotent radical of $P$.
Let $\tau^{(m)}$ denote a genuine
irreducible cuspidal representation of $GL_k^{(m)}({\A})$.
As in the previous section, construct the residue representation ${\mathcal E}_{\tau^{(m)}}^{(m)}$  defined on the group $GL_{kmc(n)}^{(m)}({\A})$. This is a Whittaker-Speh-Shalika representation of type $(km,c(n))$.
Form the Eisenstein series $E_{\tau^{(m)}}^{(m)}
(h,s)$ defined on the group $H^{(m)}({\A})$ attached to the induced
representation
\begin{align*}
\begin{cases}
\text{Ind}_{P^{(m)}({\A})}^{H^{(m)}({\A})}(
{\mathcal E}_{\tau^{(m)}}^{(m)}\otimes {\mathcal E}_{\tau^{(m)}}^{(m)}
)\delta_P^s&G=GL_n,\\
\text{Ind}_{P^{(m)}({\A})}^{H^{(m)}({\A})}{\mathcal E}_{\tau^{(m)}}^{(m)}\delta_P^s&\text{otherwise.}
\end{cases}
\end{align*}

To introduce the global integral, let $\pi_1^{(m)}$ and $\pi_2^{(m)}$ denote two genuine
 irreducible cuspidal automorphic representations defined on the group $G^{(m)}({\A})$. When $G=GL_n$, one has to be careful with the center of the group, and also with the issue of convergence. To avoid these minor complications,
in describing the integral we shall assume
 that $G$ is not the group $GL_n$.
Then the global integral we consider is
\begin{equation}\label{global1}
\int\limits_{G(F)\times G(F)\backslash G({\A})\times G({\A})}
\int\limits_{U(F)\backslash U({\A})}\varphi_1^{(m)}(g_1)\,\overline{\varphi_2^{(m)}(g_2)}\,
E_{\tau^{(m)}}^{(m)}(u(g_1,g_2),s)\,\psi_U(u)\,du\,dg_1\,dg_2.
\end{equation}
Here, $\varphi_i^{(m)}$ are vectors in the spaces of $\pi_i^{(m)}$.

We still need to check that the integral is well defined. From the description below it follows that the embeddings of the two copies of $G^{(m)}({\A})$  in $H^{(m)}({\A})$ commute with each other. It
also follows that the function
$$(g_1,g_2)\mapsto \varphi_1^{(m)}(g_1)\overline{\varphi_2^{(m)}(g_2)}
E_{\tau^{(m)}}^{(m)}((g_1,g_2),s)\ \ \ \ \ \ \ g_i\in G({\A})$$
does not depend on the choice of section $G({\A})\rightarrow G^{(m)}({\A})$. Hence the integral \eqref{global1} is well
defined and converges absolutely for $\text{Re}(s)$ large.

Finally, we describe in some detail
the unipotent group $U$, the character $\psi_U$ and the embedding of $G\times G$ inside $H$.
Let
\begin{align*}
H_{c(n)}=
\begin{cases}
GL_{2n}\times GL_{n}\times \ldots \times GL_{n}&G=GL_n,\\
Sp_{4n}&G=Sp_{2n},\\
SO_{4n}&G=SO_{2n},\\
SO_{4n+2}&G=SO_{2n+1}.
\end{cases}
\end{align*}
When $G=GL_n$, the group $GL_n$ appears $mk-1$ times in $H_{c(n)}$, and we mention this case for the sake of completeness.
Let $Q_{n,m,k}$ denote the standard parabolic subgroup of $H$ whose
Levi part is $GL_{c(n)}\times GL_{c(n)}\times\ldots GL_{c(n)}\times
H_{c(n)}$ where $GL_{c(n)}$ appears $mk-1$ times.
Let $U$ denote the unipotent radical of the parabolic subgroup
$Q_{n,m,k}$.

To define the character $\psi_U$, consider the unipotent orbit
$((2km-1)^{c(n)}1^{c(n)})$ associated with the group $H$. It follows from Collingwood
and McGovern \cite{C-M} that this is a well defined orbit for every group $H$,
and that the stabilizer of this orbit is the group $G\times G$.
From \cite{G1} we obtain that a Fourier coefficient associated with
this orbit can be constructed using the group $U$, and one can define
a character $\psi_U$ such that the stabilizer inside the Levi part
of $Q_{n,m,k}$ is the split group $G\times G$. In the next section
we shall construct this Fourier coefficient explicitly in the symplectic group case. As mentioned above, it follows
from the way the group $G\times G$ is embedded inside $H$ that the
covers in integral \eqref{global1} are compatible.

Here is our main theorem. (Notation not defined above is similar to that in other doubling method computations.)
\begin{theorem}\label{th10}
The integral \eqref{global1} is well defined, converges absolutely for $\text{Re}(s)$ large, and admits a meromorphic continuation to the whole complex plane. It is not identically zero only if
$\pi_1^{(m)}=\pi_2^{(m)}=\pi^{(m)}$. In this case, for $\text{Re}(s)$ large it is equal to
\begin{equation}\label{global2}
\int\limits_{G({\A})}\int\limits_{U_0({\A})}
<\pi^{(m)}(g)\varphi_{\pi^{(m)}},\overline{\varphi_{\pi^{(m)}}}  >f_{W({\mathcal E}_{\tau^{(m)}}^{(m)})}(\delta u_0(1,g),s)
\,\psi_U(u_0)\,du_0\,dg.
\end{equation}
Here
\begin{equation}\label{spherical}
<\pi^{(m)}(g)\varphi_{\pi^{(m)}},\varphi_{\pi^{(m)}}>=\int\limits_{G(F)\backslash G({\A})}
\varphi_{\pi^{(m)}}(g_1g)\,\varphi_{\pi^{(m)}}(g_1)\,dg_1.
\end{equation}
In particular the integral \eqref{global1} represents an Euler product.
\end{theorem}

\section{The Symplectic Group Case}\label{sp}
In this section we give some details for the symplectic group. Let $G=Sp_{2n}$.
For this group there is a  difference depending on the parity of $m$.
To give a uniform construction, let $r=m$ if $m$ is odd, and $r=m/2$ if $m$ is even. Thus, the integral we will consider is integral \eqref{global11} which agrees with integral \eqref{global1} when $m=r$
is an odd number.
Let $H=Sp_{4nrk}$. Let
${\mathcal E}_{\tau^{(r)}}^{(r)}$ be the Whittaker-Speh-Shalika representation of type $(rk,2n)$, as constructed in subsection \ref{construct}. This representation is defined on the group $GL_{2nrk}^{(r)}({\A})$, and in the notations of Definition \ref{def1} condition {\bf 1)}, we have ${\mathcal O}_{GL_{2nrk}}({\mathcal E}_{\tau^{(r)}}^{(r)})=((rk)^{2n})$.

To describe the group $U$ in integral \eqref{global11}, let $Q_{n,r,k}$ denote the parabolic subgroup of $Sp_{4nrk}$ whose Levi
part is $GL_{2n}\times\ldots\times GL_{2n}\times Sp_{4n}$. Here $GL_{2n}$ appears $rk-1$ times. Denote by $U_{n,r,k}$
or simply by $U$
the unipotent radical of $Q_{n,r,k}$. We may identify the quotient $U/[U,U]$ with the group
$$L=\text{Mat}_{2n}\oplus\ldots\oplus \text{Mat}_{2n} \oplus \text{Mat}_{2n\times 4n}.$$
Here $\text{Mat}_{2n}$ appears $rk-2$ times.
To define the character $\psi_U$ it is enough to specify it on $L$. For
$(X_1,\ldots,X_{rk-2},Y)\in L$ define $\psi_L$ as $\psi(\text{tr}(X_1+\cdots +X_{r-1})
+\text{tr}'(Y))$. Here $X_i\in \text{Mat}_{2n}$ and $Y\in \text{Mat}_{2n\times 4n}$. To define
$\text{tr}'(Y)$, write
$$Y=\begin{pmatrix} Y_1&Z_1&Y_2\\ Y_3&Z_2&Y_4 \end{pmatrix},\quad
Y_i\in \text{Mat}_{n\times n};\quad Z_j\in\text{Mat}_{n\times 2n}.$$
Then $\text{tr}'(Y)=\text{tr}(Y_1+Y_4)$.
Let $\psi_U$ denote the extension of $\psi_L$ to $U$ which is trivial on $[U,U]$.
 It follows from \cite{G1} that the
corresponding Fourier coefficient given by $U$ and $\psi_U$ is associated with the unipotent
orbit $((2rk-1)^{2n}1^{2n})$.

Finally, we specify the embedding of $(g_1,g_2)\in Sp_{2n}\times Sp_{2n}$ in $Sp_{4nrk}$. It is
given by $\text{diag}(g_1,\ldots,g_1,(g_1,g_2),g_1^*,\ldots,g_1^*)$. Here $g_1$ appears $2rk-1$ times,
and by $(g_1,g_2)$ we mean the usual embedding inside $Sp_{4n}$, i.e.\
\begin{equation}\label{diag1}
(g_1,g_2)\mapsto \begin{pmatrix} g_{1,1}&&g_{1,2}\\ &g_2&\\ g_{1,3}&&g_{1,4}\end{pmatrix},\qquad g_1=\begin{pmatrix} g_{1,1}&g_{1,2}\\ g_{1,3}&g_{1,4}\end{pmatrix}
;\quad g_{1,i}\in Mat_{n\times n}.
\end{equation}
The entries with asterisk are determined so that embedding is symplectic.
It is now easy to check that the coverings are compatible, and hence
that integral \eqref{global11} is well defined.

In \cite{G2}, Ginzburg introduced the dimension equation, which is strongly correlated with global integrals
that unfold to Euler products.
Roughly speaking, this equation states that the sum of dimensions of the representations
involved in the integral is equal to the sum of the dimensions of the groups involved.
We note that the dimension equation is satisfied in our case.
To see what is involved, let us verify this
when $m=r$ is odd.
The sum of the dimensions of the representations involved
is the sum of the dimension of the
Eisenstein series $E_{\tau^{(m)}}^{(m)}(\cdot,s)$ and the dimension of the
functional $<\pi^{(m)}(g)\varphi_{\pi^{(m)}},\varphi_{\pi^{(m)}}>$
obtained after unfolding the integral.
It follows from  \eqref{spherical} that the dimension of this functional is $\text{dim}\ Sp_{2n}$.
Thus the equation we need to verify is
\begin{equation}\label{dim1}
\text{dim}\ Sp_{2n} + \text{dim}\ E_{\tau^{(m)}}^{(m)}(\cdot,s)= 2\,\text{dim}\ Sp_{2n} +\text{dim}\ U.
\end{equation}
From \cite{G1} we see that
\begin{equation*}
\text{dim}\ E_{\tau^{(m)}}^{(m)}(\cdot,s)= \text{dim}\ {\mathcal E}_{\tau^{(m)}}^{(m)}+\text{dim}\ U(P)=\frac{1}{2}\text{dim}\ ((km)^{2n})
+ \text{dim}\ U(P).
\end{equation*}
 Here $\text{dim}\ ((km)^{2n})$ is the dimension of the unipotent orbit
$((km)^{2n})$ and $U(P)$ is the unipotent radical of the maximal parabolic $P$ defined in
Section \ref{global}. The number $\frac{1}{2}\text{dim}\ ((km)^{2n})$ is equal to
the dimension of the unipotent radical of the parabolic subgroup of $GL_{2nmk}$ whose Levi
part is $GL_{2n}\times\ldots\times GL_{2n}$ where $GL_{2n}$ appears $mk$ times. Thus its
dimension is equal to $2n^2km(km-1)$. Comparing this with the dimension of $U(P)$,
equation \eqref{dim1} follows.

We end this section with a theorem regarding the unramified computation.
To simplify notation we shall assume that $m=r$ is odd.
We first define the local $L$-functions under consideration. In general, if $\tau$ is
a local unramified representation of an $m$-fold cover of $G$, it is a constituent of an unramified principal series representation, and corresponds to some $k$-tuple $\chi=(\chi_1,\ldots,\chi_k)$, where $\chi_i$ is an unramified character of $F^*$. When $m=1$,
the values of $\chi_i$
at a local uniformizer are simply
the Satake parameters of $\tau$. For the $m$-fold cover,
we define local $L$-functions at unramified places by using the standard definition in the case of linear groups, but replacing each $\chi_i$ with its $m$-th power.

For example, let $\tau_{\nu}^{(m)}=
\text{Ind}_{B^{(m)}_{GL}}^{GL^{(m)}_k}\chi\delta_{B_{GL}}^{1/2}$ denote the local
unramified component of $\tau^{(m)}$ at a finite place $\nu$. Here $B_{GL}$ is
the Borel subgroup of $GL_k$ and $\chi=(\chi_1,\ldots,\chi_k)$ is an
unramified character. Similarly let $\pi_{\nu}^{(m)}=
\text{Ind}_{B^{(m)}_{Sp}}^{Sp^{(m)}_{2n}}\mu\delta_{B_{Sp}}^{1/2}$.
Let $p$ be a generator of the maximal ideal in the ring of integers of
the field $F_\nu$ and $q=|p|_\nu^{-1}$.
Then the local standard tensor product
$L$ function is defined by
\begin{equation}\label{tensor}
L(s,\pi_{\nu}^{(m)}\times\tau_{\nu}^{(m)})=\prod_{i=1}^{n}\prod_{j=1}^k
\frac{1}{(1-\mu_i^{m}(p)\chi_j^m(p)q^{-s})(1-\mu_i^{-m}(p)\chi_j^m(p)q^{-s})
(1-\chi_j^m(p)q^{-s})}.
\end{equation}
In a similar way one can define
the local standard  $L$-function $L(s,\tau_\nu^{(m)})$, the local exterior
square $L$-function $L(s,\tau_\nu^{(m)},\wedge^2)$ and the local symmetric square $L$-function $L(s,\tau_\nu^{(m)},\vee^2)$.

It follows from Theorem \ref{th10} that for $\text{Re}(s)$ large, integral \eqref{global2} is a product of local integrals.
At a finite unramified place $\nu$, the local integral is given by
\begin{equation}\label{local1}
Z(s,\pi^{(m)},\tau^{(m)})=\int\limits_{G}\int\limits_{U_0}
\omega_{\pi^{(m)}}(g)f_{W}(\delta u_0(1,g),s)\,\psi_U(u_0)\,du_0\,dg.\notag
\end{equation}
Here $\omega_{\pi^{(m)}}(g)$ is the spherical function of $\pi^{(m)}$, and
$f_W$ is the unramified function obtained at the place $\nu$ from the
factorizable function $f_{W({\mathcal E}_{\tau^{(m)}}^{(m)})}$. (We have dropped the subscripts $\nu$ to condense
the notation.) We have
\begin{theorem}\label{thloc1}
Suppose that $m=r$ is odd. Then  for $\text{Re}(s)$ large, the integral $Z(s,\pi^{(m)},\tau^{(m)})$ is equal to
\begin{equation}\label{local2}
\frac{L(\alpha s-\frac{\alpha-1}{2},\pi^{(m)}\times\tau^{(m)})}{L(\alpha (s-\frac{1}{2})+nm+\frac{1}{2},\tau^{(m)})\prod_{j=1}^{nm}
L(\alpha (2s-1)+2j,\tau^{(m)},\wedge^2)
\prod_{j=1}^{nm-1}
L(\alpha (2s-1)+2j+1,\tau^{(m)},\vee^2)}\notag
\end{equation}
where $\alpha=m(2nmk+1)$.
\end{theorem}

We remark that, working in the context of Brylinski-Deligne extensions, Gao \cite{Ga} also attaches Euler products to pairs of automorphic representations on covering groups.
He does so by studying the constant term of Eisenstein series on covering groups, that is,
by means of a generalization of the Langlands-Shahidi method.  However, it is not clear if his Euler products
are the same as those obtained in this note or not.

\section{Residues of Eisenstein Series on Covers of the General Linear Group}\label{res1}
In this section
we discuss Conjecture \ref{conj1} of subsection
\ref{construct} above. One way to study this Conjecture is by considering the
various constant terms of the corresponding Eisenstein series. This reduces
the problem of determining the poles to the study of the
poles of certain
intertwining operators.  This in turn reduces to the study of the poles of the
partial $L$-function $L^S(s,\tau^{(m)}\times \hat{\tau}^{(m)})$,
where $S$ is a finite set of places including all archimedean places such that
the representation $\tau^{(m)}_\nu$ is unramified for $\nu\not\in S$.
Here $\hat{\tau}^{(m)}$ is the contragredient representation of $\tau^{(m)}$.
The partial $L$-function is by definition the product over $\nu\not\in S$ of local $L$-functions
defined similarly to \eqref{tensor}.

When $m=1$, Jacquet, Piatetski-Shapiro
and Shalika studied this $L$-function by means of the Rankin-Selberg
method.  See  for example Gelbart and Shahidi \cite{G-S}, Section 1.7, for an overview of
these constructions. Using this method one may establish that
the partial $L$-function has a simple pole at $s=1$. Then it
follows from \cite{J}, for example, that Conjecture~\ref{conj1} holds.
The key property that makes this work when $m=1$ is the uniqueness of the
Whittaker model. Unfortunately, this uniqueness does not hold when $m>1$.

Suzuki studies these Eisenstein series
and their residues for higher covers in the last section of his paper \cite{Su}.
In order to establish that the residue exists, he assumes the existence of a generalized
Shimura lifting, that is, a correspondence
between irreducible cuspidal automorphic representations of the group
$GL_n^{(m)}({\A})$ and automorphic representations of the group
$GL_n({\A})$ which satisfies certain properties.
At the moment, such a lifting has only been proved in full
for covers of $GL_2$ (Flicker, \cite{F}).
If one had the lift in general and if the lifted automorphic representation of $GL_n({\A})$
was cuspidal, then it would follow from the case $m=1$  that the partial $L$-function
$L^S(s,\tau^{(m)}\times \hat{\tau}^{(m)})$ has a simple pole at $s=1$. This would then
imply Conjecture \ref{conj1}. We summarize:
\begin{proposition}\label{prop10}
Assume that $\pi^{(m)}$ satisfies the two conditions in \cite{Su} ,  Section 8.5, p. 752. Then Conjecture~\ref{conj1} holds.
\end{proposition}

An alternative method to study the residues of Eisenstein is given in
Jacquet and Rallis \cite{J-R}. In
this method one constructs a certain global nonzero integral which involves the residue representation. The non-vanishing of the global
period then implies that the residue is nonzero. When $m=2$ we propose
the following construction. For simplicity let us consider the case of
the maximal parabolic Eisenstein series.
Let $E_{\tau^{(2)}}^{(2)}(g,s)$ denote the Eisenstein series defined on the group $GL_{2k}
^{(2)}({\A})$, which is associated with the induced representation
\begin{equation*}
\text{Ind}_{P_{2,k}^{(2)}({\A})}^{GL_{2k}^{(2)}({\A})}
(\tau^{(2)}\otimes \tau^{(2)})\delta_{P_{2,k}}^{s}.\notag
\end{equation*}
 Let $E_{\tau^{(2)}}
^{(2)}(g)$ denote the residue of this series at the point
$s=(2k+1)/4k$. Our goal is to prove that the representation generated by these residues is nonzero.  Let
$\Theta^{(2)}_{Sp}$  denote the theta representation of the group $Sp_{2k}^{(2)}({\A})$.
It is the minimal representation of this group. We have the following result.
\begin{proposition}\label{period1}
Let $\theta_{Sp}^{(2)}$ be a vector in the space of $\Theta^{(2)}_{Sp}$.  Then
the integral
\begin{equation}\label{per1}
\int\limits_{Sp_{2k}(F)\backslash Sp_{2k}({\A})}E_{\tau^{(2)}}
^{(2)}(g)\,\theta_{Sp}^{(2)}(g)\,dg\notag
\end{equation}
converges absolutely.  Moreover, for some choice of data the integral is not identically zero.
In particular the residue representation is nonzero.
\end{proposition}
Unfortunately, at this point we do not have a way to extend this result
to higher order covers.

\end{document}